\documentclass{ifacconf}

\usepackage{graphicx}      % include this line if your document contains figures
\usepackage{natbib}        % required for bibliography
\usepackage{booktabs}
\usepackage{amssymb}
\usepackage{amsmath,mathtools}
\usepackage{amsfonts}
\usepackage{bbm}
\usepackage{aligned-overset}
\usepackage{caption}
\usepackage[dvipsnames]{xcolor}

\begin{document}
\begin{frontmatter}

\title{Distributionally robust model predictive control for wind farms}
\thanks[footnoteinfo]{\copyright 2023 the authors. This work has been accepted to IFAC for publication under a Creative Commons Licence CC-BY-NC-ND.}

\author[First]{Christoph Mark} 
\author[First]{Steven Liu} 

\address[First]{Institute  of  Control  Systems,  Department  of  Electrical  and  Computer
	Engineering, University of Kaiserslautern-Landau, 67663
	Kaiserslautern, Germany (e-mail: \{cmark, steven.liu\}@rptu.de).}

\begin{abstract}                % Abstract of not more than 250 words.
In this paper, we develop a distributionally robust model predictive control framework for the control of wind farms with the goal of power tracking and mechanical stress reduction of the individual wind turbines. We introduce an ARMA model to predict the turbulent wind speed, where we merely assume that the residuals are sub-Gaussian noise with statistics contained in an moment-based ambiguity set. We employ a recently developed distributionally model predictive control scheme to ensure constraint satisfaction and recursive feasibility of the control algorithm. The effectiveness of the approach is demonstrated on a practical example of five wind turbines in a row.
\end{abstract}

\begin{keyword}
Predictive control, Constrained control, Stochastic control
\end{keyword}

\end{frontmatter}
%===============================================================================
\section{Introduction}
A large part of green energy production is currently covered by wind farms (WF), where several wind turbines (WT) are placed in close proximity to each other to reduce the cost of cabling and maintenance. One problem that occurs in such an environment is that each wind turbine generates a wake that moves downstream and is characterized by a flow velocity deficit and increased turbulence intensity \cite{barthelmie2007modelling}. The flow velocity deficit directly impacts the power production of downstream turbines \cite{barthelmie2010quantifying}, while the increased turbulence intensity increases the fatigue loads \cite{bossuyt2017measurement}. 

In this paper, a distributionally robust model predictive controller (DR-MPC) is developed as a supervisory controller for a wind farm with the primary objective of dynamically allocating a required wind farm power reference to the individual wind turbines in the field. The WT power references are then tracked by underlying local WT controllers, which operate on a much faster timescale (millisecond range) compared to the WF controller (second range). A secondary objective of the wind farm controller is to reduce fatigue loads of the turbines to increase their overall lifetime. The DR-MPC algorithm is based on our previous publication \cite{mark2022recursively}, but has been extended to include cost functions for output variables.

\textit{Related work:}
In \cite{riverso2016model}, the authors consider the same setup as we do and use a stochastic MPC (SMPC) to design a supervisory control system for wind farms, adopting the probabilistic SMPC framework from \cite{farina2013probabilistic}. However, their approach is based on the assumption that the true wind speed is normally distributed and the moments are known exactly. In \cite{boersma2019stochastic}, a scenario-based SMPC for power reference tracking is developed, where Gaussianity of the wind speed distribution is assumed. Fatigue load reduction is not considered explicitly in this work. The authors of \cite{spudic2011wind} investigated a deterministic MPC approach for wind farm control. Similar to our approach, the goal was to track power and reduce mechanical stress, however, the stochasticity of the wind is neglected and assumed to be constant over the prediction horizon. This approach was extended to a distributed MPC in \cite{spudic2015cooperative}. In terms of wind turbine control, several papers have been published that address fatigue reduction, such as \cite{evans2014robust}, where a robust MPC was developed for oscillation damping, or \cite{gros2017real}, where an economic nonlinear MPC was applied to reduce structural and actuator fatigue.
\subsection{Notation}
A probability space is defined by the triplet $(\Omega, \mathcal{F}, \mathbb{P})$, where $\Omega$ is the sample space, $\mathcal{F}$ the Borel $\sigma$-algebra on $\Omega$ and $\mathbb{P}$ the probability measure on $(\Omega, \mathcal{F})$.  
The set of all probability distributions supported on $\Omega$ with finite second moment is $\mathcal{M}(\Omega)$. Given an event $E_1$ we define the probability occurrence as $\mathbb{P}(E_1)$ and the conditional probability given $E_2$ as $\mathbb{P}(E_1 | E_2)$. For a random variable $w$, we define the expected value as $\mathbb{E}(w)$, whereas the conditional expectation of $w$ conditional to a random variable $x$ is denoted as $\mathbb{E}(w | x)$. The weighted 2-norm w.r.t.~a positive definite matrix $Q = Q^\top$ is $\Vert x \Vert_Q^2 = x^\top Q x$. Positive definite and semidefinite matrices are indicated as $A\succ0$ and $A\succeq0$, respectively. The pseudo inverse of a matrix $A$ is denoted as $A^\dagger$. The stacked column vector $x \in \mathbb{R}^{N n}$ of subvectors $x_1, \ldots, x_N \in \mathbb{R}^{n}$ is defined as $x = \mathrm{col}_{i = 1, \ldots, N}(x_i)$. The Kroneker product is denoted as~$\otimes$.
\subsection{Outline}
In Section \ref{sec:problem}, we introduce the wind farm model and pose the general optimization problem of interest. Section \ref{sec:power_system:MPC_formulation} is devoted to the theoretical background of the distributionally robust MPC, which is based on our previous publication \cite{mark2022recursively}. In Section \ref{sec:power_system:simulations}, we perform two simulations using a wind farm with five wind turbines in a row. The paper closes with Section \ref{sec:conclusion}, where we summarize the results and provide a brief outlook.
\section{Problem description}
\label{sec:problem}
In this paper, we use a linearized version of the NREL WT as proposed by \cite{riverso2016model}, where the $i$-th WT is described by a linear time-invariant system of the form
\begin{subequations}
	\label{eq:power_system:wind_turbine_ss}
	\begin{align}
		  x_i(k+1) &= A_i x_i(k) + B_i u_i(k) + E_i w_i(k) \label{eq:power_system:WT:state_equation} \\
		  y_i(k) &= C_i x_i(k) + D_i u_i(k) + F_i w_i(k) \label{eq:power_system:WT:out_equation}
	\end{align}
\end{subequations}
with state $x_i = (\beta_i, \omega_{i,\mathrm{r}}, \omega^\mathrm{f}_{i, \mathrm{g}}) - (\beta_{i,0}, \omega_{i, \mathrm{r0}}, \omega_{i,\mathrm{g0}})$, input $u_i = P^{\mathrm{wt}}_{i, \mathrm{ref}} - P^{\mathrm{wt}}_{i, 0}$, disturbance $w_i = \tilde{w}_i - w_{i,0}$ and output $y_i = (F_{i,\mathrm{t}}, T_{i, \mathrm{s}}) - (F_{i,\mathrm{t0}}, T_{i, \mathrm{s0}})$. The individual components are the blade pitch angle $\beta_i \: [\mathrm{^\circ}]$, the rotor angular velocity $\omega_{i,\mathrm{r}} \: [\mathrm{rad/s}]$, the filtered generator angular velocity $\omega_{i, \mathrm{g}} \: [\mathrm{rad/s}]$, the power reference $P^{\mathrm{wt}}_{i, \mathrm{ref}} \: [\mathrm{W}]$, the effective wind velocity $\tilde{w}_i \: [\mathrm{m/s}]$, the tower bending force $F_{i,\mathrm{t}} \: [\mathrm{N}]$ and the main shaft torque $T_{i,\mathrm{s}} \: [\mathrm{Nm}]$. Quantities with a subindex $0$ denote the nominal operating (linearization) point. The wind farm model is obtained by stacking $N_{\mathrm{wt}} \in \mathbb{N}$ individual WT models \eqref{eq:power_system:wind_turbine_ss}, so that $x = \mathrm{col}_{i = 1, \ldots, N_{\mathrm{wt}}}(x_i)$, $u = \mathrm{col}_{i = 1, \ldots, N_{\mathrm{wt}}}(u_i)$, $w = \mathrm{col}_{i = 1, \ldots, N_{\mathrm{wt}}}(w_i)$ and $y = \mathrm{col}_{i = 1, \ldots, N_{\mathrm{wt}}}(y_i)$, while the dynamic matrices $A - F$ are given by block diagonal stacking of the local matrices $A_i - F_i$ for all $i = 1, \ldots, N_{\mathrm{wt}}$, resulting in
\begin{subequations}
	\label{eq:power_system:wind_farm_ss}
	\begin{align}
		x(k+1) &= A x(k) + B u(k) + E w(k) \label{eq:power_system:WF:state_equation} \\
		y(k) &= C x(k) + D u(k) + F w(k). \label{eq:power_system:WF:out_equation}
	\end{align}
\end{subequations}
The main objective of a wind farm controller in the above rated region is to distribute the wind power reference provided by the system operator to each wind turbine in the field while minimizing fatigue load \cite{knudsen2015survey, andersson2021wind}. Fatigue loads result from repetitive stress reversals on a specific part of the structure, where typical fatigue prone components are the turbine tower and the generator shaft \cite{spudic2012coordinated}. Therefore, we formulate the following infinite horizon stochastic optimal control problem
\begin{subequations}
	\label{eq:smpc_optimization}
	\begin{align}
		\!\min_{u(k) \forall k \in \mathbb{N}} & \quad \mathbb{E}_{\mu^*} \left( \sum_{k=0}^{\infty} \Vert y(k) \Vert_{Q_\mathrm{y}}^2 + \Vert u(k) \Vert_{R}^2  \bigg | x(0) \right)  \label{eq:smpc:cost} \\
		\mathrm{s.t.} \hspace{0.91em} & \quad x(k+1) = A x(k) + B u(k) + E w(k) \nonumber \\
		& \quad  y(k) = C x(k) + D u(k) + F  w(k) \\
		& \quad \mathbb{P}(l_{j}^\top u(k) \leq 1 \: | \: x(0)) \geq p^u_{j} \quad j \in \{1, \ldots, s\} \label{eq:smpc:input_chance_constraints} \\
		& \quad  \mathbbm{1}^\top  u(k) = 0, \label{eq:smpc:input} 
	\end{align}
\end{subequations}
where \eqref{eq:smpc:cost} is an expected value quadratic cost function that penalizes deviations of the output and input with weights $Q_\mathrm{y} \succeq 0$ and $R \succ 0$, \eqref{eq:smpc:input_chance_constraints} denotes a set of individual input chance constraints of probability level $p_j^{\mathrm{u}} \in (0,1)$ and \eqref{eq:smpc:input} enforces that the power deviations in sum are equal to zero to cover the nominal demand.

 The optimization problem \eqref{eq:smpc_optimization} contains several sources of intractability, i.e., (i) the control input $u$ in the presence of an additive uncertainty renders the problem infinite dimensional; (ii) the expectation in \eqref{eq:smpc:cost} is taken w.r.t.~the unknown probability distribution $\mu^*$ and (iii) the chance constraints \eqref{eq:smpc:input_chance_constraints} are evaluated under the unknown probability measure $\mathbb{P}$. Point (i) will be tackled in Section~\ref{sec:SADF}, where we introduce a simplified affine disturbance feedback (SADF) parameterization, while the uncertainty sources (ii) and (iii) are addressed with a distributionally robust cost function and distributionally robust constraints in Sections \ref{sec:cost_fcn} and \ref{sec:chance_constraints} based on a moment-based ambiguity set introduced in Section \ref{sec:ambiguity}.

\section{Distributionally robust MPC}
\label{sec:power_system:MPC_formulation}
In this section, we use a DR-MPC scheme recently proposed by the authors \cite{mark2022recursively} to approximate the infinite horizon problem \eqref{eq:smpc_optimization} over a finite prediction horizon.

\subsection{Prediction dynamics}
To distinguish between closed-loop and predicted states and inputs, we introduce the $N$-step ahead prediction of \eqref{eq:power_system:WF:state_equation} over a horizon of length $N \in \mathbb{N}$
\begin{align}
	\bar{x}_k = \bar{A} x_{0|k} + \bar{B} \bar{u}_k + \bar{E} \bar{w}_k, \label{eq:state_prediction}
\end{align}
where $\bar{x}_k = [x^\top_{0|k}, x^\top_{1|k}, \ldots, x^\top_{N|k}]^\top$ denotes the state sequence, $\bar{u}_k = [u^\top_{0|k}, \ldots, u^\top_{N-1|k}]^\top$ the input sequence and $\bar{w}_k = [w^\top_{0|k}, \ldots, w^\top_{N-1|k}]^\top$ the disturbance sequence, while the matrices are defined as
{\scriptsize
	\begin{align*}
		&\bar{A} = \begin{bmatrix}
			I \\
			A \\
			\vdots \\
			A^N
		\end{bmatrix}\hspace{-0.3em}, 
		\bar{B} = \begin{bmatrix}
			0 & 0 & \dots & 0 \\
			B & 0 & \dots & 0 \\
			AB & B & \dots & 0 \\
			\vdots & \ddots & \ddots & 0 \\
			A^{N-1}B & \dots & AB & B
		\end{bmatrix}\hspace{-0.3em},
		\bar{E} = \begin{bmatrix}
			0 & 0 & \dots & 0 \\
			E & 0 & \dots & 0 \\
			AE & E & \dots & 0 \\
			\vdots & \ddots & \ddots & 0 \\
			A^{N-1}E & \dots & AE & E
		\end{bmatrix}.
	\end{align*}
}%
Similar to the state sequence \eqref{eq:state_prediction}, we represent the output equation \eqref{eq:power_system:WF:out_equation} in a compact form as
\begin{align}
	\bar{y}_k = \bar{C}  \bar{x}_k + \bar{D}  \bar{u}_k + \bar{F}  \bar{w}_k, \label{eq:output_prediction}
\end{align}
where $ \bar{y}_k = [  y_{0|k}^\top, \ldots, y_{N|k}^\top]^\top$ and the matrices are given by
{\small
	\begin{align*}
		&\bar{C} = \begin{bmatrix}
			C& 0 & \dots  & 0& 0 \\
			0 & C & \dots & 0 & 0 \\
			\vdots & \vdots & \ddots & \vdots & \vdots \\
			0 & 0 & \dots & C & 0\\
			0 & 0 & \dots & 0 & C
		\end{bmatrix}, \bar{D} = \begin{bmatrix}
			0 & 0 & \dots & 0 \\
			D& 0 & \dots & 0 \\
			0 & D & \dots & 0 \\
			\vdots & \vdots & \ddots & \vdots \\
			0 & 0 & \dots & D 
		\end{bmatrix}, 
	\bar{F} = \begin{bmatrix}
			0 & 0 & \dots & 0 \\
			F& 0 & \dots & 0 \\
			0 & F & \dots & 0 \\
			\vdots & \vdots & \ddots & \vdots \\
			0 & 0 & \dots & F 
		\end{bmatrix}.
	\end{align*}
}%
Note that $\bar{w}_k$ represents a turbulent wind speed prediction that is typically non-i.i.d.~and correlated in time. Therefore, to restore the sub-Gaussianity of the random variables as required by \cite{mark2022recursively}, we identify an auto-regressive moving average model that serves as a whitening filter for the turbulent wind speed.

\subsection{ARMA model and ambiguity set} 
\label{sec:ambiguity}
An ARMA model represents a stochastic process in terms of two polynomials, where the first one represents the auto-regressive (AR) part and the second one the moving average (MA) part \cite{box2015time}. In particular, an $\mathrm{ARMA}(p,q)$ model with $p$ AR terms and $q$ MA terms is given by
\begin{align*}
	 {w}(k) = \sum_{l= 1}^p a_l  {w}(k-l) + \sum_{l=1}^q b_l \epsilon(k-l) + \epsilon(k),
\end{align*}
where $\epsilon$ is a zero-mean i.i.d~white noise. In related work, e.g. \cite{ono2013risk, riverso2016model}, the authors make a more stringent assumption that the noise $\epsilon$ is normally distributed, which in case of wind turbulence data is prone to be wrong, cf. \cite{van2015distributionally}. Therefore, we treat $\epsilon$ as a zero-mean white noise with unknown (but finite) variance $\Sigma_\epsilon \succ 0$. In practice, one needs to identify the ARMA model with limited data.~Therefore, the empirical variance is typically falsified due to sample inaccuracies, for which we introduce a moment-based ambiguity set that captures the true variance with high probability
\begin{align}
	\label{eq:power_system:ambiguity_set}
	\mathcal{P}(w_{\mathrm{0}}, T_\mathrm{I}) \coloneqq \left\{  \mu \in \mathcal{M}(\mathbb{R}^n) \  \middle\vert \begin{array}{l}
		\mathbb{E}_{\mu}( \epsilon )  = 0 \\
		\mathbb{E}_{\mu}( \epsilon \epsilon^\top ) \preceq \kappa^{(w_{\mathrm{0}}, T_\mathrm{I})}_\beta \hat{\Sigma}_\epsilon^{(w_{\mathrm{0}}, T_\mathrm{I})}
	\end{array}\right\}.
\end{align}
Note that we parameterize the ambiguity set with the mean wind speed and turbulence intensity pair $(w_{\mathrm{0}}, T_\mathrm{I})$. The ambiguity radius $\kappa^{(w_{\mathrm{0}}, T_\mathrm{I})}_\beta$ can readily be found with \cite[Prop. 1]{mark2022recursively}.

We identify for each wind turbine $i = 1, \ldots, N_\mathrm{wt}$ an $\mathrm{ARMA}(p, p-1)$ model, which can be converted to a canonical form similar to \cite{ono2013risk}, i.e,
\begin{align*}
	\psi_i(k+1) &= A_{\psi, i} \psi_i(k) + B_{\psi,i} \epsilon_i(k)  \\
	 {w}_i(k) &= C_{\psi,i} \psi_i(k), 
\end{align*}
where the matrices are defined as follows
\begin{align*}
	A_{\psi,i} \coloneqq \begin{bmatrix}
		a_{i,1} & 1 & 0 & \dots & 0 \\
		a_{i,2} & 0 & 1 &       & 0 \\
		\vdots & \vdots & & \ddots &  \\
		a_{i,p-1} & 0 & 0 & & 1 \\
		a_{i,p} & 0 & 0 & \dots & 0
	\end{bmatrix}, 
	B_{\psi,i} \coloneqq \begin{bmatrix}
		1 \\ b_{i,1} \\ \vdots \\ b_{i,p-2} \\ b_{i,p-1}
	\end{bmatrix},
\end{align*}
$C_{\psi,i} \coloneqq \begin{bmatrix}
	1 & 0 & \cdots & 0
\end{bmatrix}$ and the auxiliary state vector $\psi_i$ is given by $\psi_i(k) = \begin{bmatrix}
	  w^\top_i(k), \psi^\top_{i,2}(k), \dots, \psi^\top_{i,p}(k)
\end{bmatrix}^\top$ with 
\begin{align*}
	\psi_{i,j}(k) &= \sum_{l = j}^p a_{i,l}  {w}_i(k+j-l-1) \\
	&+ \sum_{l = j-1}^{p-1} b_{i,l} \epsilon_i(k+j-l-1) \quad \forall i \in \{1, \ldots, N_\mathrm{wt}\}.
\end{align*}
To obtain farm-wide wind predictions, we stack the local matrices and vectors, such that $A_\psi = \mathrm{diag}(A_{\psi,1}, \ldots, A_{\psi,N_\mathrm{wt}})$, $B_\psi = \mathrm{diag}(B_{\psi,1}, \ldots, B_{\psi,N_\mathrm{wt}})$, $C_\psi = \mathrm{diag}(C_{\psi,1}, \ldots, C_{\psi,N_\mathrm{wt}})$ and $\psi = \text{col}_{i = 1}^{N_\mathrm{wt}}(\psi_i)$. A $N$-step prediction of the turbulent wind speed is readily given by
\begin{align}
	  \bar{w}_k \coloneqq \bar{C}_\psi \bar{A}_\psi \psi(k) + \bar{C}_\psi\bar{B}_\psi \bar{\epsilon}_k, \label{eq:power_syste:arma_prediction}
\end{align}
where $\bar{C}_\psi \coloneqq \mathrm{diag}(C_\psi, \ldots, C_\psi)$,
{\small
\begin{align*}
	&\bar{A}_\psi \coloneqq \begin{bmatrix}
		I \\
		A_\psi \\
		A_\psi^2 \\
		\vdots \\
		A_\psi^{N-1}
	\end{bmatrix}, \bar{B}_\psi \coloneqq \begin{bmatrix}
		0 & 0 & \dots & 0 \\
		B_\psi & 0 & \dots & 0 \\
		A_\psi B_\psi & B_\psi & \dots & 0 \\
		\vdots & \ddots & \ddots & 0 \\
		A_\psi^{N-2}B_\psi & \dots & A_\psi B_\psi  & B_\psi
	\end{bmatrix},
\end{align*}
}%
while the random vector $\bar{\epsilon}_k$ is zero-mean and each element is i.i.d.~with variance $\Sigma_\epsilon^{(w_{\mathrm{0}}, T_\mathrm{I})}$. 
%\begin{rem}
%	As a byproduct of the ARMA model, a covariance reduction of the new random variable $\epsilon$ compared to the original random variable $w$ is usually achieved. This aspect is important for the DR-MPC implementation because the chance constraints and the cost function depend directly on the covariance matrix, i.e., the lower the covariance, the lower the conservatism.
%\end{rem}

\subsection{Simplified affine disturbance feedback} 
\label{sec:SADF}
To render the resulting MPC optimization problem finite dimensional, we parameterize the control input with a SADF policy, cf. \cite{zhang2020stochastic}, of the form 
\begin{align}
	\bar{u}_k = \bar{v}_k + \bar{M}_k \bar{\epsilon}_k, \label{eq:power_system:disturbance_feedback}
\end{align}
where the matrices are defined as
\begin{align*}
	&\bar{M}_k \coloneqq \begin{bmatrix}
		0 & 0 & \dots & 0 \\
		M_{1|k} & 0 & \dots & 0 \\
		\vdots & \ddots & \ddots & 0 \\
		M_{N-1|k} & \dots & M_{1|k} & 0
	\end{bmatrix}, \:
	\bar{v}_k \coloneqq \begin{bmatrix}
		v_{0|k} \\
		v_{1|k} \\
		\vdots \\
		v_{N-1|k}
	\end{bmatrix}.
\end{align*}
Next, we substitute the state prediction \eqref{eq:state_prediction}, the ARMA prediction \eqref{eq:power_syste:arma_prediction} and the SADF policy \eqref{eq:power_system:disturbance_feedback} into the output prediction \eqref{eq:output_prediction}, resulting in
\begin{align}
	\bar{y}_k &= \bar{C} \bar{x}_k + \bar{D} \bar{u}_k + \bar{F} \bar{w}_k \nonumber \\
	\overset{\eqref{eq:state_prediction}, \eqref{eq:power_syste:arma_prediction}, \eqref{eq:power_system:disturbance_feedback}}&{=} \underbrace{\bar{C} \bar{A} x_{0|k} + (\bar{C} \bar{B} + \bar{D}) \bar{v}_k + (\bar{C} \bar{E} + \bar{F}) \bar{C}_\psi \bar{A}_\psi \psi(k)}_{  \tilde{\bar{y}}_k} \nonumber\\
	& \quad + \underbrace{[\bar{C} \bar{B} \bar{M}_k  + \bar{D} \bar{M}_k + (\bar{C} \bar{E} + \bar{F}) \bar{C}_\psi\bar{B}_\psi ]}_{\Psi_k} \bar{\epsilon}_k. \label{eq:power_system:output_prediction}
\end{align}
The output prediction now depends only on the initial conditions $x_{0|k}, \psi(k)$, as well as on the optimization variables $\bar{v}_k$ and $\bar{M}_k$.
\begin{rem}
	The decision variables using the SADF policy grow linearly in the prediction horizon $N$, whereas the original affine disturbance feedback policy grows quadratically \cite{zhang2020stochastic}. Thus, the SADF policy results in less demanding optimization problems.
\end{rem}
\subsection{Cost function}
\label{sec:cost_fcn}
We approximate the infinite horizon cost function \eqref{eq:smpc:cost} over the prediction horizon $N$, while we replace the expected value over $\mu^*$ with the supremum over all distributions contained in the ambiguity set \eqref{eq:power_system:ambiguity_set}, i.e.,
\begin{align}
	J_k&= \sup_{\mu \in \mathcal{P}}\: \mathbb{E}_{\mu} \left(  \bar{y}_k^\top \bar{Q}_\mathrm{y}  \bar{y}_k +  \bar{u}_k^\top \bar{R}  \bar{u}_k \bigg | x(k)\right) \nonumber\\
	\overset{\eqref{eq:power_system:disturbance_feedback}, \eqref{eq:power_system:output_prediction}}&{=}  \text{tr}\left(\hat{\Sigma}_N^{(w_{\mathrm{0}}, T_\mathrm{I})} \big[\bar{H}_{\mathrm{y},k}^\top \bar{Q}_\mathrm{y} \bar{H}_{\mathrm{y},k} + \bar{H}_{\mathrm{u},k}^\top \bar{R} \bar{H}_{\mathrm{u},k} \big]  \right), \label{eq:power_system:cost_reformulation}
\end{align}
where $\bar{H}_{\mathrm{y},k} = \begin{bmatrix}
	\Psi_k^\top &  \tilde{\bar{y}}_k^\top
\end{bmatrix}^\top$, $\bar{H}_{\mathrm{u},k} = \begin{bmatrix}
	\bar{M}_k^\top &  \bar{v}_k^\top
\end{bmatrix}^\top$, $\bar{Q}_\mathrm{y} = \mathrm{diag}(Q_\mathrm{y}, \ldots, Q_\mathrm{y})$ and $\bar{R}= \mathrm{diag}(R, \ldots, R)$. The worst-case covariance matrix $\hat{\Sigma}_N^{(w_{\mathrm{0}}, T_\mathrm{I})}$ is defined through the moment-based ambiguity set \eqref{eq:power_system:ambiguity_set}
{\normalsize
\begin{align*}
	\hat{\Sigma}_N^{(w_{\mathrm{0}}, T_\mathrm{I})} &\coloneqq \sup_{\mu \in \mathcal{P}} \bigg( \begin{bmatrix}
		\bar{\epsilon}_k \\ 1
	\end{bmatrix} \begin{bmatrix}
		\bar{\epsilon}_k \\ 1
	\end{bmatrix}^\top \bigg| x(k) \bigg) \\
	\overset{\mathrm{i.i.d.}}&{=}  
	\mathrm{diag}\left(I_N \otimes \sup_{\mu \in \mathcal{P}(w_{\mathrm{0}}, T_\mathrm{I})} \bigg( \begin{bmatrix}
		\epsilon \\ 1
	\end{bmatrix} \begin{bmatrix}
		\epsilon\\ 1
	\end{bmatrix}^\top \bigg| x(k) \bigg), 1\right) \\
	\overset{\eqref{eq:power_system:ambiguity_set}}&{=} \mathrm{diag}(I_N \otimes \kappa^{(w_{\mathrm{0}}, T_\mathrm{I})}_\beta \hat{\Sigma}_\epsilon^{(w_{\mathrm{0}}, T_\mathrm{I})}, 1),
\end{align*}
}%
where the first equality follows from the i.i.d.~sequence $\bar{\epsilon}_k$, i.e., $\bar{\epsilon}_k$ contains $N$-times the i.i.d.~random variable $\epsilon$. For details on the reformulation steps, please refer to our recent paper \cite{mark2022recursively}.
\subsection{Chance constraints}
\label{sec:chance_constraints}
Since the probability measure $\mathbb{P}$ required for the chance constraints \eqref{eq:smpc:input_chance_constraints} is unknown, we instead impose distributionally robust chance constraints, i.e., we enforce the chance constraint for all distributions contained in the ambiguity set \eqref{eq:power_system:ambiguity_set}, resulting in
	\begin{align*}
		&\inf_{\mu \in \mathcal{P}} \:  \mathbb{P}(l_{t,j}^\top \bar{u}_k \leq 1 \: | \: x(k)) \geq p^u_{s}	\nonumber\\
		&\overset{\eqref{eq:power_system:disturbance_feedback}}{=}\inf_{\mu \in \mathcal{P}} \: \mathbb{P} \bigg( \begin{bmatrix}
			\bar{\epsilon}_k^\top & 1
		\end{bmatrix}
		\begin{bmatrix}
				l^\top_{t,j} \bar{M}_k \\
				l^\top_{t,j} \bar{v}_k
		\end{bmatrix} \leq 1 \bigg \vert x(k)\bigg) \geq p_j^\mathrm{u},
	\end{align*}
which can equivalently be expressed as second-order cone constraint via \cite[Thm 3.1]{calafiore2006distributionally}
	\begin{align}
		l^\top_{t,j} \bar{v}_k \leq 1 - \sqrt{\frac{p_u}{1 - p_u}} \Vert l_{t,j}^\top \bar{M}_k (\hat{\Sigma}_N^{(w_{\mathrm{0}}, T_\mathrm{I})})^{\frac{1}{2}} \Vert_2.
	\end{align}
	The vector $l_{t,j}$ is a lifted version $l_j$ from constraint \eqref{eq:smpc:input_chance_constraints} to fit the dimension of the nominal input vector $\bar{v}_k$.
\subsection{MPC optimization problem}
At each time step $k\in \mathbb{N}$, we solve the following MPC optimization problem
{\small
\begin{subequations}
	\label{eq:mpc_problem}
	\begin{align}
		&\!\min_{\bar{v}_k, \bar{M}_k, \lambda_k}  && \text{tr}\left(\hat{\Sigma}_N^{(w_{\mathrm{0}}, T_\mathrm{I})} \big[\bar{H}_{\mathrm{y},k}^\top \bar{Q}_\mathrm{y} \bar{H}_{\mathrm{y},k} + \bar{H}_{\mathrm{u},k}^\top \bar{R} \bar{H}_{\mathrm{u},k} \big]  \right) \label{eq:power_system:mpc:cost} \\
		&\quad\mathrm{s.t.} & & \bar{y}_k = \tilde{\bar{y}}_k + \Psi_k \bar{\epsilon}_k\\
		& & & \bar{z}_k = \bar{A} x_{0|k} + \bar{B} \bar{v}_k + \bar{E} \bar{C}_\psi \bar{A}_\psi \psi(k) \label{eq:mpc:state}\\
		& & &  x_{0|k} = (1-\lambda_k) x(k) + \lambda   z^*_{1|k-1}, \quad \lambda_k \in [0, 1] \label{eq:power_system:mpc:init_constraints}	\\
		& & & l^\top_{t,j} \bar{v}_k \leq 1 - \sqrt{\frac{p_u}{1 - p_u}} \Vert l_{t,j}^\top \bar{M}_k (\hat{\Sigma}_N^{(w_{\mathrm{0}}, T_\mathrm{I})})^{\frac{1}{2}} \Vert_2 \label{eq:power_system:mpc:input_constraints}\\
				& & & \mathbbm{1}^\top v(k) = 0 \quad \forall j \in \{1, \ldots, s\} \: \forall t \in \{0, \ldots, N - 1\}, \label{eq:power_constraint}
\end{align}
\end{subequations}
}%
where \eqref{eq:power_system:mpc:init_constraints} is a so-called interpolating initial constraint that ensures recursive feasibility if the optimization problem is feasible at time $k = 0$, cf. \cite{kohler2022recursively, mark2022recursively}. The constraint \eqref{eq:mpc:state} denotes the nominal state prediction, while \eqref{eq:power_constraint} enforces that the mean power demand is in sum equal to zero.

The optimal solution of the MPC optimization problem \eqref{eq:mpc_problem} yields the SADF pair $( \bar{v}_k^*, \bar{M}_k^*)$ and the mean state prediction $\bar{z}_k^*$. Following the lines of \cite{mark2022recursively}, we obtain an equivalent admissible error feedback control policy via
\begin{align*}
	\bar{K}_k^* &= (I + \bar{M}_k^* \bar{E}^\dagger \bar{B})^{-1} \bar{M}^*_k \bar{E}^\dagger \\
	 \bar{g}_k^* &= (I + \bar{M}_k^* \bar{E}^\dagger \bar{B})^{-1} (  \bar{v}^*_k - \bar{M}_k^* \bar{E}^\dagger A   z^*_{0|k}),
\end{align*}
where 
\begin{align*}
	{\tiny
		\bar{K}_k \coloneqq \begin{bmatrix}
			K_{0|k} & 0 & \ldots & 0 & 0 \\
			K_{1|k} & K_{0|k} & \ldots & 0 & 0 \\
			\vdots & \ddots & \ddots & \vdots & 0 \\
			K_{N-1|k} & \ldots & K_{1|k} & K_{0|k} & 0 
		\end{bmatrix}, \bar{g}_k \coloneqq \begin{bmatrix}
			g_{0|k} \\
			g_{1|k} \\
			\vdots \\
			g_{N-1|k}
		\end{bmatrix},
	}
\end{align*}
while the input to the wind turbines is defined as
\begin{align*}
	P^{\mathrm{wt}}_{\mathrm{ref}}(k) = u(k) = P^{\mathrm{wt}}_{\mathrm{ref0}} +   g^*_{0|k} + K_{0|k}^* (  x(k) -   z_{0|k}^*),
\end{align*}
where $P^{\mathrm{wt}}_{\mathrm{ref0}}$ is a vector of nominal power references for all wind turbines.
\section{Numerical example}
\label{sec:power_system:simulations}
\begin{figure}[!b]
	\centering
	\includegraphics[width=0.5\linewidth]{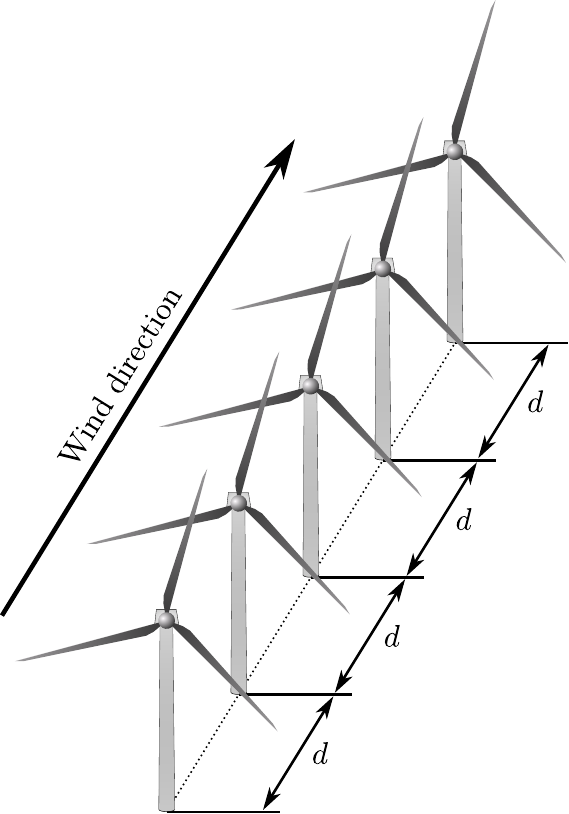}
	\caption{Simulation setup with five equidistantly arranged wind turbines.}
	\label{fig:power_system:simulation_setup}
\end{figure}

In the following, we apply our proposed DR-MPC to a wind farm consisting of $N_\mathrm{wt} = 5$ NREL $5$ MW wind turbines in a row, see Figure \ref{fig:power_system:simulation_setup}, where each WT is equidistantly arranged with $d = 400$ m. We use the Matlab/Simulink toolbox SimWindFarm (SWF) \cite{grunnet2010aeolus} as our simulation environment.

\paragraph*{Controller types} We compare the DR-MPC to an open-loop \textit{Scheduler} that assigns a constant power references $P^{\mathrm{wt}}_{i,\mathrm{ref}} = 3 \: \mathrm{MW}$ for each wind turbine $i = 1, \ldots, N_\mathrm{wt}$ for the entire simulation horizon of $T = 900$ s, i.e., the wind farm should nominally produce $P^{\mathrm{wf}}_{\mathrm{ref0}} = 15 \: \mathrm{MW}$. In addition, we consider the \textit{SWF controller} \cite{grunnet2010aeolus}, which dynamically distributes the power references based on the available wind power estimates of each wind turbine
\begin{align*}
	P_{i, \mathrm{avail}} = \min\left\{ P_{\mathrm{0}}, \frac{1}{2} \pi \rho R^2 \tilde{w}_i^3 C^{\mathrm{max}}_{\mathrm{p}}\right\} \quad \forall i \in \{1, \ldots, N_\mathrm{wt}\},
\end{align*}
where $\tilde{w}_i$ is the measured (effective) wind speed at turbine $i$, $P_\mathrm{0}$ the rated power and $C^{\mathrm{max}}_{\mathrm{p}} = 0.45$ the maximum power coefficient. Therefore, the SWF controller distributes the power as follows
\begin{align*}
	P^{\mathrm{wt}}_{i, \mathrm{ref}} = \max\left\{0, \min \bigg\{ P_{\mathrm{0}}, \frac{P^{\mathrm{wf}}_{\mathrm{ref0}} P_{i, \mathrm{avail}}}{\sum_{i = 1}^{N_\mathrm{wt}} P_{i, \mathrm{avail}}} \bigg\} \right\}
\end{align*}
for all $i = 1, \ldots, N_\mathrm{wt}$.

\paragraph*{Performance metrics} Similar to \cite{riverso2016model}, we use the following metrics to evaluate the performance of each controller:
\begin{itemize}
	\item Tracking $\displaystyle J_\mathrm{p} = \sqrt{ \frac{1}{T} \sum_{k = 0}^{T-1} \sum_{i = 1}^{N_{\mathrm{wt}}} \frac{(P_{i, \mathrm{out}}(k) - P^{\mathrm{wt}}_{i,\mathrm{ref}}(k))^2}{N_{\mathrm{wt}} P_{\mathrm{0}}} }$,
	\item Shaft fatigue $\displaystyle J_\mathrm{s} = \mathrm{std}\left(\frac{\sum_{k = 0}^{T-1} \sum_{i = 1}^{N_{\mathrm{wt}}}T_{i,\mathrm{s}}(k)}{N_{\mathrm{wt}}T_{\mathrm{s0}}} \right)$,
	\item Tower fatigue $\displaystyle J_\mathrm{t} = \mathrm{std}\left(\frac{\sum_{k = 0}^{T-1} \sum_{i = 1}^{N_{\mathrm{wt}}} F_{i,\mathrm{t}}(k)}{N_{\mathrm{wt}} T_{\mathrm{t0}}} \right)$,
\end{itemize}
where $P_{i, \mathrm{out}}$, $T_{i,\mathrm{s}}$ and $F_{i,\mathrm{t}}$ denote the power output, main shaft torque and tower bending force of WT $i$, while $T_{\mathrm{s0}} =  2.5 \cdot 10^6$ and $T_{\mathrm{t0}} = 0.27 \cdot 10^6$ are standardization constants. To reduce the tuning effort of the MPC cost function, we fix the output weight $\bar{Q}_\mathrm{y}$ to
{\small
\begin{align*}
	\bar{Q}_\mathrm{y} = I_{N+1} \otimes \mathrm{diag} \left( \begin{bmatrix}
	 \frac{1}{F_{\mathrm{t0}}^2 N} & 0 \\
		0 &  \frac{100}{T_\mathrm{s0}^2 N}
	\end{bmatrix}, \ldots, \begin{bmatrix}
		 \frac{1}{F_{\mathrm{t0}}^2 N} & 0 \\
		0 &  \frac{100}{T_\mathrm{s0}^2 N}
	\end{bmatrix} \right)
\end{align*}
}%
and, analogously, the input weighting matrix to
\begin{align*}
	\bar{R} = I_N \otimes \begin{bmatrix}
		\frac{r}{P^2_{\mathrm{0}} N} & \cdots & 0 \\
		\vdots & \ddots & \vdots \\
		0 & \cdots & \frac{r}{P^2_{\mathrm{0}} N}
	\end{bmatrix},
\end{align*}
where $r \in \mathbb{R}_{> 0}$. Thus, it remains to tune the parameter $r$, which introduces a trade-off between tracking performance and fatigue load reduction. We consider a prediction horizon of $N=5$ seconds for each simulation.
\paragraph*{Operation scenario}
We consider a realistic wind farm operation scenario in which some wind turbines temporarily operate in the below rated region due to deficiencies in wind speed, while others operate in the above rated region, cf. \cite{morgan2011probability}. The wind field has a mean velocity of $w_{\mathrm{0}}~=~12~\mathrm{m/s}$ and a turbulence intensity of $T_\mathrm{I} = 0.1$, yielding a turbulence variance of $\sigma^2_{\mathrm{  w}} = 1.44$.
For the DR-MPC, we identify during the offline phase for each wind turbine an ARMA model based on an independent wind scenario (training data) of $1000$ time steps, i.e., $N_s = 1000$ samples. We identified both an ARMA$(2,1)$ and an ARMA$(3,2)$ model and chose the latter due to its lower root mean square error between the ARMA prediction and the training data.

In view of \cite[Prop. 1]{mark2022recursively}, we derive an ambiguity radius of $\kappa^{(w_{\mathrm{0}}, T_\mathrm{I})}_\beta = 2.36$ with a confidence of $1-\beta = 0.95$, while the empirical covariance matrix of the ARMA residuals is given by
\begin{align*}
	\hat{\Sigma}_\epsilon^{(12, 0.1)} = \mathrm{diag}(0.255, \:  0.270, \: 0.288, \: 0.262, \:0.274).
\end{align*}
We constrain the input deviations to $\pm 1 \: \mathrm{MW}$ around the nominal operating point of $3 \: \mathrm{MW}$ with a probability of $90 \%$, which allows the DR-MPC to dynamically dispatch the power references depending on the available wind speed, while ensuring a power tracking goal and minimizing fatigue load. This is enforced with the input chance constraint \eqref{eq:power_system:mpc:input_constraints}, while in addition a penalty term  $5\lambda^2$ is added to the cost function \eqref{eq:power_system:mpc:cost}, which enforces that the interpolated initial constraint \eqref{eq:power_system:mpc:init_constraints} favors the feedback initialization.

\subsection{Simulation results}
In Table \ref{power_system:tab:scenario2}, we compare the performance of each controller regarding tracking accuracy, tower fatigue and transmission shaft fatigue, where the Scheduler is considered to be the baseline, i.e., $100 \%$. Numbers below $100 \%$ indicate a relative increase in performance, while numbers above $100 \%$ reflect a relative decrease in performance. In particular, for $r=1$, we increase the tracking performance compared to the scheduler by approximately $52.6 \%$ and compared to the SWF controller by $4.4 \%$. The tracking performance increase comes at the price of increasing the tower fatigue by $27.8 \%$, while reducing the main shaft fatigue by $10 \%$. A reasonable choice is $r=500$, which only marginally increases the mechanical stress on the tower, while still increasing the tracking performance by nearly~$34 \%$. 

\begin{figure}[!t]
	\centering
	\includegraphics[width=1\linewidth]{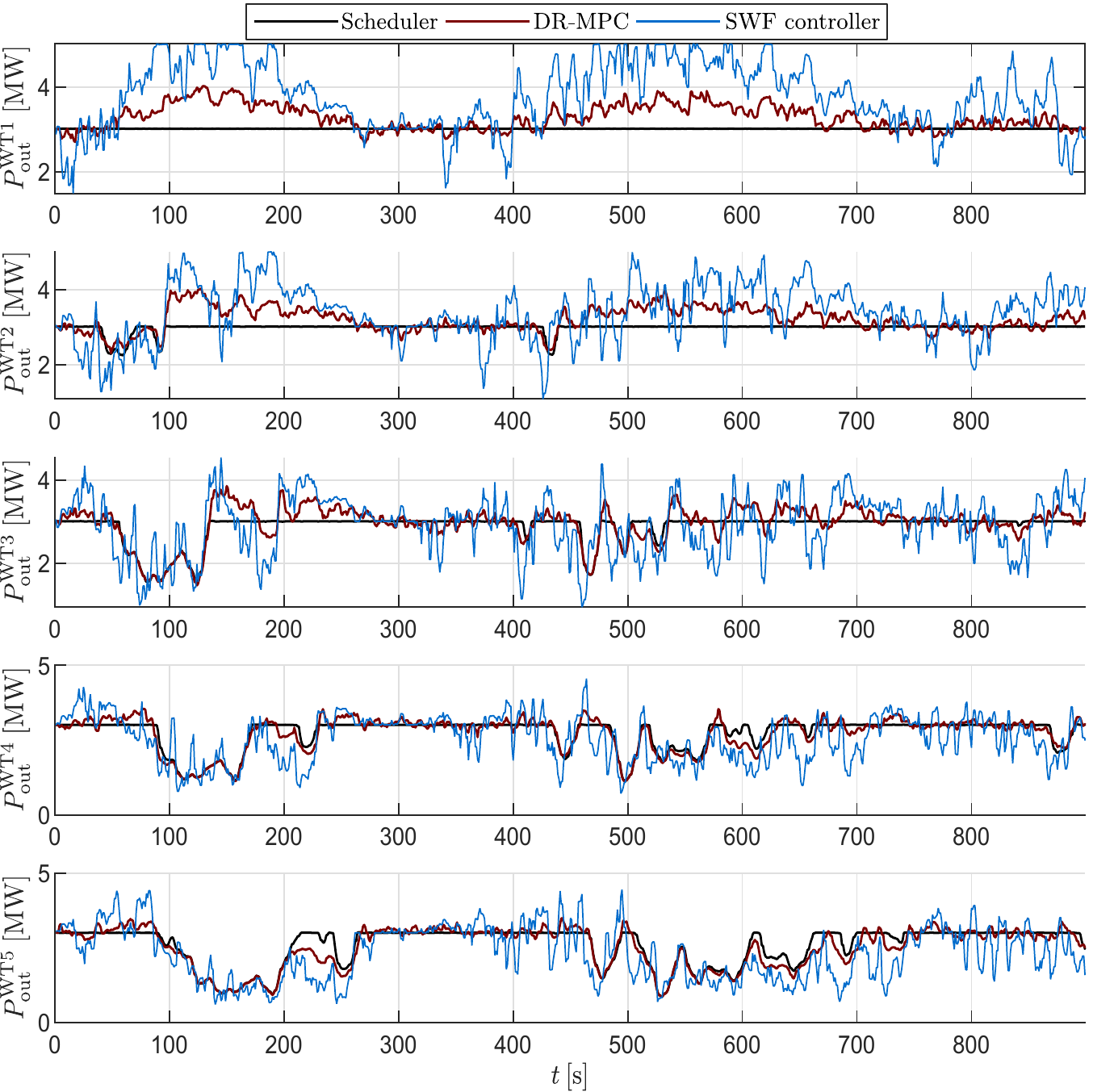}
	\caption{Wind turbine electrical power output. The DR-MPC uses the weight $r=500$.}
	\label{fig:power_system:simulation:scenario21:power_output}
\end{figure}

In Figure \ref{fig:power_system:simulation:scenario21:power_output}, we illustrate the electrical power output of each wind turbine. First, a wake turbulence can be observed through the wind farm that causes the power output of the downstream turbines to temporarily fall below the nominal value of $3 \mathrm{MW}$, e.g., at time $t~\approx~50$ the second turbine experiences a wind speed deficit that causes the power output to fall below the nominal value, while at time $t \approx 95$ the wake turbulence affects turbine five. This illustrates the weakness of the open-loop scheduler, as there is no immediate feedback to increase the output of the upstream turbines and compensate for the loss of the others. The SWF controller and the DR-MPC are both feedback strategies, which dynamically allocate power references, i.e., turbine one increases its power output as soon as the downstream turbines drop below their nominal value. This immediately results in an increase in tracking performance, but also additional mechanical fatigue. Compared to the SWF controller, the DR-MPC allows for systematically balancing the tracking performance with the increasing mechanical fatigue, which can be seen in less variance in the electrical output among all turbines. 

The DR-MPC optimization problem \eqref{eq:mpc_problem} is implemented with Yalmip \cite{lofberg2004yalmip} and Mosek \cite{aps2019mosek}, and is solved in $35.6$ ms on average on an Intel i7-9700k processor with 16gb ram, confirming real-time capability as the controller should typically operate in the seconds range \cite{spudic2015cooperative}.

\begin{center}
	\captionof{table}{Performance comparison for different controller types/parameterizations.}
	\begin{tabular}{@{}lrrr@{}}
		\toprule     
		Method    			& $J_\mathrm{p}$ 		   		& $J_\mathrm{t}$  	& $ J_\mathrm{s} $ \\
		\midrule
		Scheduler  			& $0.0999$    	&  $0.3217$   	& $0.0734$\\
		SWF controller  	& $51.79 \%$ 	&  $131.99 \%$ 	& $90.15 \%$ \\ 
		DR-MPC $R = 1$      & $47.38 \%$ 	&  $127.85 \%$    & $89.99 \%$ \\
		DR-MPC $R = 500$   & $65.79 \%$ 	&  $105.30 \%$   & $93.41 \%$ \\
		DR-MPC $R = 10^3$   & $84.07 \%$ 	&  $101.88 \%$    & $93.67 \%$ \\
		DR-MPC $R = 10^4$   & $98.55 \%$ & $100.12 \%$ &  $99.21 \%$  \\
		\bottomrule                          
	\end{tabular} 
	\label{power_system:tab:scenario2}
\end{center}

\section{Conclusion}
\label{sec:conclusion}
In this paper, we presented a distributionally robust MPC approach to tackle the problem of coordinating individual wind turbines inside of a wind farm. The main objective hereby was to ensure power tracking, while a secondary goal was to reduce the mechanical stress acting on the tower and main transmission shaft. 
In a case study of five wind turbines in a row, we numerically verified the increase in tracking performance as well as the reduction in mechanical fatigue compared to a simple open-loop scheduler approach. In addition, we illustrated the trade-off between power tracking and fatigue reduction. 
We considered an ARMA model to predict the turbulent wind speed locally for each wind turbine individually, neglecting the broader picture of spatial correlations of the wind field. This can be improved by considering a spatio-temporal wind speed forecast that includes wind measurements from neighboring turbines, e.g., as proposed by \cite{zhu2018wind}. This could further increase the tracking performance, i.e., when a wind deficit is measured at an upstream turbine, it is inevitably passed on to the downstream turbines, allowing us to anticipate the temporary drop in output power in the future. Therefore, the output of the unaffected wind turbines can be increased to compensate for the loss of output power of the others.

\bibliography{ifacconf}             % bib file to produce the bibliography

\end{document}